\documentclass[a4paper,10pt]{amsart}
\usepackage[german,english]{babel}
\usepackage[latin1]{inputenc}
\usepackage{ucs}
\usepackage{amsmath}
\usepackage{amsfonts}
\usepackage{amssymb}
\usepackage{amsthm}
\usepackage{fontenc}
\usepackage{graphicx}
\usepackage{url}
\usepackage{color}
\usepackage[T1]{fontenc}
\usepackage{a4wide}
\usepackage{xcolor}     
\usepackage{mdframed}   
\newenvironment{graybox}%
  {\begin{mdframed}[backgroundcolor=lightgray,linecolor=lightgray]}%
  {\end{mdframed}}
\newtheorem{definition}{Definition}
\newtheorem{remark}{Remark}

\title[Optimal heating of an indoor swimming pool]{
Optimal heating of an indoor swimming pool
}

\author{Monika Wolfmayr}

\address[Monika Wolfmayr]{
Faculty of Information Technology, University of Jyväskylä,
P.O. Box 35 (Agora), FIN-40014 Jyväskylä, Finland}
\email{monika.k.wolfmayr@jyu.fi}

\begin{document}

\begin{abstract}This work presents the derivation of a model for the heating process of the air
of a glass dome, where an indoor swimming pool is located in the bottom of the dome. 
The problem can be reduced from
a three dimensional to a two dimensional one. 
The main goal is the formulation of a
proper optimization problem for computing the optimal heating of the air after a given time.
For that, the model of the heating process as 
a partial differential equation is formulated 
as well as the optimization problem subject to the time-dependent partial differential equation. This yields
the optimal heating of the air under the glass dome such that the desired temperature distribution is 
attained after a given time. 
The discrete formulation of the optimization problem 
and a proper numerical method for it,  
the projected gradient method, are discussed. 
Finally, numerical experiments are presented which show the practical
performance of the optimal control problem and its numerical solution method discussed.
\end{abstract}

\maketitle

\section{Introduction}
\label{mw:sec:1:introduction}

Modeling the heating of an object is an important task in many applicational problems.
Moreover, a matter of particular interest is to find the optimal heating of an object such that it has a 
desired temperature distribution after some given time. 
In order to formulate such optimal control problems and to solve them,
a cost functional subject to a time-dependent partial differential equation (PDE) is derived.
One of the profound works paving the way for PDE-constrained optimization's relevance
in research and application during the last couple of decades
is definitely Lion's work  \cite{MW:Lions:1971} from 1971.
Some recent published monographs discussing PDE-constrained optimization as well as various efficient
computational methods for solving them are, e.g., \cite{MW:BorziSchulz:2012},
\cite{MW:HinzePinnauUlbrichUlbrich:2009}, and \cite{MW:Troeltzsch:2010}, where the latter one
is used as basis for the discussion on solving the optimal heating problem of this work.

The goal of this work is to derive a simple mathematical model for finding the optimal heating of the air in
a glass dome
represented by a half sphere, where a swimming pool is located in the bottom of the dome 
and the heat sources (or heaters) 
are situated on a part of the boundary of the glass dome.
The process from the model to the final numerical simulations usually involves several steps. The main steps in this
work are the setting up of the mathematical model for the physical problem, 
obtaining some analytical results
of the problem, presenting a proper discretization for the continuous problem and finally computing the 
numerical solution of the problem.
The parabolic optimal control problem is discretized by the finite element method in space,
and in time, we use the implicit Euler method for performing the time stepping.
The used solution algorithm for the discretized problem is the projected gradient method,
which is for instance applied in \cite{MW:HerzogKunisch:2010} as well as in more detail discussed in
\cite{MW:GruverSachs:1980, MW:HinzePinnauUlbrichUlbrich:2009, MW:Kelley:1999, MW:NocedalWright:1999}.

We want to emphasize that the model and the presented optimization methods
for the heating process of this work are only one example for a possible modeling and solution. 
In fact, the stated model problem has potential
for many modeling tasks for students and researchers. 
For instance, different material parameters for the dome as well
as for air and water could be studied more carefully.
The optimal modeling of the heat sources could 
be stated as a shape optimization problem or instead of optimizing the temperature of the air
in the glass dome, one could optimize the water temperature, which would correspond
to a final desired temperature distribution correponding to the boundary of the glass dome, where
the swimming pool is located, for the optimal control problem.
Another task for the students could be to compute many simulations with, e.g.,
Matlab's \texttt{pdeModeler} to derive a better understanding of the problem in the pre-phase of studying 
the problem of this work. However, we only want to mention here a few other possibilities for
modeling, studying and solving the optimal heating problem amongst many
other tasks, and we are not focusing on 
them in the work presented here.

This article is organized as follows:
First, the model of the heating process is formulated in Section \ref{mw:sec:2:model}.
Next, Section \ref{mw:sec:3:OCP} introduces 
the optimal control problem, which
describes the optimal heating of the glass dome 
such that the desired temperature distribution is attained after a given time. 
In Section \ref{mw:sec:4:exun}, proper function spaces are presented in order to 
discuss existence and uniqueness of the optimal control problem there as well.
We derive the reduced optimization problem in Section \ref{mw:sec:5:reduced-problem}
before discussing its discretization and the numerical method for solving it,
the projected gradient method, 
in more detail
in Section \ref{mw:sec:6:num-method}.
Numerical results 
are presented 
as well as conclusions are drawn in the final Section \ref{mw:sec:7:num-results}.

\section{Modeling}
\label{mw:sec:2:model}

This section presents the modeling process.
The physical problem is described in terms of mathematical language,
which includes
formulating
an initial version 
of the problem, but then simplifying it in order
to derive a version of the problem which is easier to solve. However, at the same time,
the problem has to be 
kept accurate enough in order
to compute an approximate solution being close enough to the original solution.
That is exactly one of the major goals of mathematical modeling.
In the following, we introduce the domain describing the glass dome, where an
indoor swimming pool is located in the bottom of the dome,
and the position of the heaters.
The 
concrete
equations describing the process of heating and the cost functional
subject to them 
modeling the optimization task are discussed in the next Section \ref{mw:sec:3:OCP}.

We have an indoor swimming pool which is located under a glass dome. For simplicity, it is assumed
that we have an isolated system in the glass dome, so no heat can leak from the domain.
The swimming pool covers the floor of the glass dome and
we assume that the heaters are placed next to the floor up on the glass all around the dome.
The target of the minimization functional is to reach a desired temperature distribution at the
end of a given time interval $(0,T)$, where $T > 0$ denotes the final time, with the least possible cost. 

\begin{figure}[bht!]
\includegraphics[scale=.51]{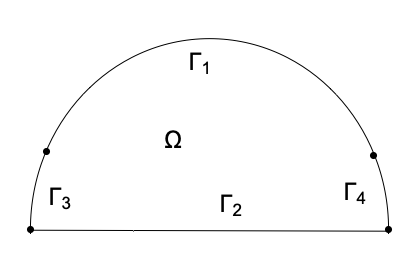}
\caption{The domain $\Omega$ reduced from 3d to 2d due to symmetry properties
describing 
the glass dome and its heaters placed at the ground of the glass
dome next to the floor, hence subdividing the boundary $\Gamma = \partial \Omega$
 into four parts
$\Gamma_1$, $\Gamma_2$, $\Gamma_3$ and $\Gamma_4$.}
\label{mw:fig:1:modeldomain2d}    
\end{figure}

Due to the symmetry properties of the geometry as well as the uniform distribution of the water temperature,
we reduce the three dimensional (3d) problem to a two dimensional (2d) one.
The dimension reduction makes the numerical computations more simple.
In the following, the 2d domain is denoted by $\Omega$ and its boundary by
$\Gamma = \partial \Omega$. We assume that
$\Omega \subset \mathbb{R}^2$ is a bounded
Lipschitz domain.
We subdivide its boundary $\Gamma$ into four parts:
the glass part $\Gamma_1$, the floor which is the swimming pool $\Gamma_2$ and the
heaters $\Gamma_3$ and $\Gamma_4$.
The domain $\Omega$ and its boundaries are illustrated in Figure \ref{mw:fig:1:modeldomain2d}.

\section{Optimal control problem}
\label{mw:sec:3:OCP}

In this section, the 
optimal control problem 
is formulated, where an optimal
control function $u$ has to be obtained corresponding to the heating of the heat sources on the
boundary $\Gamma_R := \Gamma_3 \cup \Gamma_4$
such that the state $y$ reaches a desired temperature distribution $y_d$ after a given
time $T$. This problem can be formulated in terms of a PDE-constrained optimization problem,
which means minimizing a cost functional subject to a PDE
and with $u$ being the control function.

Let $Q_T := \Omega \times (0,T)$ denote the space-time cylinder
with the lateral surface $\Sigma := \Gamma \times (0,T),$ where $T > 0$ denotes
the final time.
The optimal control problem is given as follows:
\begin{equation}
\label{mw:eq:costfunc}
\min_{(y,u)} J(y,u) = \frac{1}{2} \int_{\Omega} (y(\boldsymbol{x},T) - y_d(\boldsymbol{x}))^2 \, d\boldsymbol{x}
+ \frac{\lambda}{2} \int_0^T \int_{\Gamma_R} u(x,t)^2 \, ds \, dt
\end{equation}
such that
\begin{align}
y_t - \triangle y &= 0 &&\quad \text{in } Q_T := \Omega \times (0,T),
\label{mw:eq:PDE1} \\
\frac{\partial y}{\partial n} &= 0 &&\quad \text{on } \Sigma_1 := \Gamma_1 \times (0,T),
\label{mw:eq:PDE2} \\
y &= g &&\quad \text{on } \Sigma_2 := \Gamma_2 \times (0,T),
\label{mw:eq:PDE3} \\
\frac{\partial y}{\partial n} + \alpha y &= \beta u &&\quad \text{on } \Sigma_R := \Gamma_R \times (0,T),
\label{mw:eq:PDE4} \\
y(0) &= y_0 &&\quad \text{in } \Omega,
\label{mw:eq:PDE5}
\end{align}
where $g$ is the given constant water temperature,
$\lambda \geq 0$ is the cost coefficient or control parameter, and
$\alpha$ and $\beta$ are constants describing the heat transfer, which are
modeling parameters and have to be chosen carefully.
The control $u$ denotes the radiator heating, which has to be chosen within a certain temperature range. Hence,
we choose the control functions from the following set of admissible controls:
\begin{equation}
\label{mw:eq:admissible}
u \in U_{\text{ad}} = \{v \in L^2(\Sigma_R) : u_a(\boldsymbol{x},t) \leq u(\boldsymbol{x},t) \leq u_b(\boldsymbol{x},t)
\text{ a.e. on } \Sigma_R\},
\end{equation}
which means that $u$ has to fulfill so called box constraints.
The equations \eqref{mw:eq:PDE2}, \eqref{mw:eq:PDE3} and \eqref{mw:eq:PDE4}
are called Neumann, Dirichlet and Robin boundary conditions, respectively.
Equation \eqref{mw:eq:PDE2} characterises a no-flux condition in normal direction.
Equation \eqref{mw:eq:PDE3} describes a constant temperature distribution. 
Regarding equation \eqref{mw:eq:PDE4}, $\alpha = \beta$ would be a reasonable choice from
the physical point of view because this would mean 
that the temperature increase 
at this part of the boundary is proportional to the difference between the temperature there and outside.
However, a decoupling of the parameters makes sense too,
see \cite{MW:Troeltzsch:2010}, and does not change anything for the actual discussion and computations,
since $\alpha = \beta$ could be chosen at any point.

The goal is to find the optimal set of state and control $(y,u)$ such that the cost is minimal.

\section{Existence and uniqueness}
\label{mw:sec:4:exun}

In this section, we discuss some basic results on the existence and uniqueness of the
parabolic initial-boundary value problem \eqref{mw:eq:PDE1}--\eqref{mw:eq:PDE5}, whereas
we exclude the details. They can be found in \cite{MW:Troeltzsch:2010}.
We first introduce proper function spaces leading 
to a setting, where existence and uniqueness of the solution can be proved.
\begin{definition}
The normed space $W^{1,0}_2(Q_T)$ is defined as follows 
\begin{equation}
W^{1,0}_2(Q_T) = \{ y \in L^2(Q_T) : D_i y \in L^2(Q_T) \, \forall i = 1, ..., d  \}
\end{equation}
with the norm
\begin{equation}
\|y\|_{W^{1,0}_2(Q_T)} = \left( \int_0^T \int_{\Omega} (|y(\boldsymbol{x},t)|^2 
+ |\nabla y(\boldsymbol{x},t)|^2 ) \, d\boldsymbol{x} \, dt \right)^{1/2},
\end{equation}
where $D_i y$ denotes the spatial derivative of $y$ in $i$-direction and $d$ is the
spatial dimension.
\end{definition}
For the model problem of this work, the dimension is $d=2$.
In the following, let $\{V, \|\cdot\|_V\}$ be a real Banach space.
More precisely, we will consider $V = H^1(\Omega)$ in this work.
\begin{definition}
The space $L^p(0,T;V)$,
$1 \leq p < \infty$, denotes the linear space of all equivalence classes of measurable
vector valued functions $y: [0,T] \rightarrow V$ such that 
\begin{equation}
\int_0^T \|y(t)\|_V^p \, dt < \infty.
\end{equation}
The space $L^p(0,T;V)$ is a Banach space with respect to the norm
\begin{equation}
\|y\|_{L^p(0,T;V)} := \left(\int_0^T \|y(t)\|_V^p \, dt \right)^{1/p}.
\end{equation}
\end{definition}
\begin{definition}
The space $W(0,T) = \{y \in  L^2(0,T;V): y' \in L^2(0,T;V^*)\}$ 
is equipped with the norm
\begin{equation}
\|y\|_{W(0,T)} = \left( \int_0^T (|y(t)|_V^2 
+ |y'(t)|_{V^*}^2 ) \, dt \right)^{1/2}.
\end{equation}
It is a Hilbert space with the scalar product
\begin{equation}
(u,w)_{W(0,T)} =  \int_0^T (u(t),w(t))_V \, dt + \int_0^T (u'(t),w'(t))_{V^*} \, dt.
\end{equation}
\end{definition}
%
The relation $V \subset H = H^* \subset V^*$ is called a
Gelfand or evolution triple and describes a chain of dense and continuous embeddings.


The problem \eqref{mw:eq:PDE1}--\eqref{mw:eq:PDE5}
has a unique weak solution $y \in W^{1,0}_2(Q_T)$ for a given $u \in U_{\text{ad}}$. 
Moreover, the solution depends continuously on the data, which means that
there exists a constant $c > 0$ being independent of $u$, $g$ and $y_0$ such that
\begin{equation}
\label{mw:estimate:W102}
\max_{t \in [0,T]} \|y(\cdot,t)\|_{L^2(\Omega)} + \|y\|_{W^{1,0}_2(Q_T)}
\leq c (\|u\|_{L^2(\Sigma_R)} + \|g\|_{L^2(\Sigma_2)} + \|y_0\|_{L^2(\Omega)})
\end{equation}
for all $u \in L^2(\Sigma_R)$, $g \in L^2(\Sigma_2)$ and $y_0 \in L^2(\Omega)$.
Hence, problem \eqref{mw:eq:PDE1}--\eqref{mw:eq:PDE5}
is well-posed in $W^{1,0}_2(Q_T)$.
Furthermore, since $y \in W^{1,0}_2(Q_T)$ and it is a weak solution of
problem \eqref{mw:eq:PDE1}--\eqref{mw:eq:PDE5},
$y$ also belongs to $W(0,T)$.
The following estimate holds:
\begin{equation}
\label{mw:estimate:W0T}
\|y\|_{W(0,T)} 
\leq \tilde c (\|u\|_{L^2(\Sigma_R)} + \|g\|_{L^2(\Sigma_2)} + \|y_0\|_{L^2(\Omega)})
\end{equation}
for some constant $\tilde c > 0$ being independent of $u$, $g$ and $y_0$.
Hence, problem \eqref{mw:eq:PDE1}--\eqref{mw:eq:PDE5}
is also well-posed in the space $W(0,T)$.

Note that the Neumann boundary conditions on $\Gamma_1$ are included in both
estimates
\eqref{mw:estimate:W102} and \eqref{mw:estimate:W0T}
related to the well-posedness of the problem
(as discussed in  \cite{MW:Troeltzsch:2010}).
However, the Neumann boundary conditions \eqref{mw:eq:PDE2} are equal to zero.

Under the assumptions that
$\Omega \subset \mathbb{R}^2$ is a bounded Lipschitz domain with boundary $\Gamma$,
$\lambda \geq 0$ is a fixed constant, $y_d \in L^2(Q_T)$, $\alpha, \beta \in L^{\infty}(\Sigma_R)$, and
$u_a, u_b \in L^2(\Sigma_R)$ with $u_a \leq u_b$ a.e. on $\Sigma_R$,
together with the existence and uniqueness result on the
parabolic initial-boundary value problem \eqref{mw:eq:PDE1}--\eqref{mw:eq:PDE5}
in $W(0,T)$, the optimal control problem \eqref{mw:eq:costfunc}--\eqref{mw:eq:admissible}
has at least one optimal control $\bar u \in U_{\text{ad}}$. In case of $\lambda > 0$ the optimal control
$\bar u$ is uniquely determined.

\section{Reduced optimization problem}
\label{mw:sec:5:reduced-problem}

In order to solve the optimal control problem
\eqref{mw:eq:costfunc}--\eqref{mw:eq:admissible}, we derive the
so called reduced optimization problem first.

Since the problem is well-posed as discussed in the previous section, we can formally eliminate the state
equation \eqref{mw:eq:PDE1}--\eqref{mw:eq:PDE5}
and the minimization problem reads as follows
\begin{equation}
\label{mw:eq:costfunc-reduced}
\min_{u} \bar{J}(u) = \frac{1}{2} \int_{\Omega} (y_u(\boldsymbol{x},T) - y_d(\boldsymbol{x}))^2 \, d\boldsymbol{x}
+ \frac{\lambda}{2} \int_0^T \int_{\Gamma_R} u(x,t)^2 \, ds \, dt
\end{equation}
such that \eqref{mw:eq:admissible} is satisfied. The problem \eqref{mw:eq:costfunc-reduced} is called
reduced optimization problem. Formally, the function $y_u$ denotes that the state function is depending on $u$.
However, for simplicity we can set again $y_u = y$.
In order to solve the problem \eqref{mw:eq:costfunc-reduced}, we apply the projected gradient method.
The gradient of $\bar{J}$ has to be calculated by deriving the adjoint problem which is given by 
\begin{align}
\label{mw:eq:PDE1red} 
-p_t - \triangle p &= 0 &&\quad \text{in } Q_T, \\
\label{mw:eq:PDE2red}
\frac{\partial p}{\partial n} &= 0 &&\quad \text{on } \Sigma_1, \\
\label{mw:eq:PDE3red}
p &= 0 &&\quad \text{on } \Sigma_2, \\
\label{mw:eq:PDE4red} 
\frac{\partial p}{\partial n} + \alpha p &= 0 &&\quad \text{on } \Sigma_R, \\
\label{mw:eq:PDE5red}
p(T) &= y(T) - y_d &&\quad \text{in } \Omega.
\end{align}
The gradient of $\bar{J}$ is given by
\begin{equation}
\label{mw:eq:gradJbar}
\nabla \bar{J}(u(\boldsymbol{x},t)) = \beta \chi_{\Gamma_R} p(\boldsymbol{x},T-t) + \lambda u(\boldsymbol{x},t)
\end{equation}
with $\chi_{\Gamma_R}$ denoting the characteristic function on $\Gamma_R$.
The projected gradient method can be now applied 
for computing the solution of the PDE-constrained optimization problem \eqref{mw:eq:costfunc}--\eqref{mw:eq:admissible}.
We denote by
\begin{equation}
\label{mw:eq:projection-admissible}
\mathcal{P}_{[u_a,u_b]}(u) = \max\{u_a, \min\{u_b,u\}\}
\end{equation}
the projection onto the set of admissible controls $U_{\text{ad}}$.

Now putting everything together, the optimality system for \eqref{mw:eq:costfunc}--\eqref{mw:eq:admissible}
and a given $\lambda > 0$ 
reads as follows \\
\begin{graybox}
\begin{equation}
\label{mw:eq:OS}
\begin{aligned}
y_t - \triangle y = 0 &\qquad \qquad -p_t - \triangle p = 0 \quad \qquad \,\,\,\,\, \text{in } Q_T \\
\frac{\partial y}{\partial n} = 0  &\qquad \qquad \frac{\partial p}{\partial n} = 0 \qquad \qquad \quad \quad \,\,\,\, \text{on } \Sigma_1 \\
y = g &\qquad \qquad p = 0 \qquad \qquad \qquad \quad \,\,\, \text{on } \Sigma_2 \\
\frac{\partial y}{\partial n} + \alpha y = \beta u &\qquad \qquad
\frac{\partial p}{\partial n} + \alpha p = 0 \qquad \qquad \, \, \, \text{on } \Sigma_R \\
y(0) = y_0 &\qquad \qquad p(T) = y(T) - y_d \qquad \quad \text{in } \Omega, \\
u = &\, \mathcal{P}_{[u_a,u_b]}(-\frac{1}{\lambda} \beta p).
\end{aligned}
\end{equation}
\end{graybox}

\vspace{0.5cm}
In case that $\lambda = 0$, the projection formula changes to
\begin{equation}
\label{mw:eq:proj-lambda0}
\begin{aligned}
u(\boldsymbol{x},t) &= u_a(\boldsymbol{x},t), \qquad \text{ if } \beta(\boldsymbol{x},t) p(\boldsymbol{x},t) > 0, \\
u(\boldsymbol{x},t) &= u_b(\boldsymbol{x},t), \qquad \text{ if } \beta(\boldsymbol{x},t) p(\boldsymbol{x},t) < 0.
\end{aligned}
\end{equation}

\begin{remark}
In case that there are no control constraints imposed, the projection formula simplifies to
$u = - \lambda^{-1} \beta p$.
\end{remark}

\section{Discretization and numerical method}
\label{mw:sec:6:num-method}

In order to numerically solve the optimal control problem
\eqref{mw:eq:costfunc}--\eqref{mw:eq:admissible}, which is equivalent
to solving \eqref{mw:eq:OS},
we discretize the heat equation in space
by the 
finite element method and in time, we use the implicit Euler method
for performing the time stepping.

We approximate the functions $y, u$ and $p$ by finite element functions
$y_h, u_h$ and $p_h$ from the conforming finite element
space $ V_h = \mbox{span} \{\varphi_1, \dots, \varphi_n\} $ with the basis functions
$\{\varphi_i(\boldsymbol{x}): i=1,2,\dots,n_h \}$, where
$h$ denotes the discretization parameter with $n = n_h = \mbox{dim} \, V_h = O(h^{-2})$.
We use standard, continuous, piecewise linear finite elements and a regular triangulation $\mathcal{T}_h$ 
to construct the finite element space $V_h$. 
For more information, we refer the reader to
\cite{MW:Ciarlet:1978} as well as to the newer publications
\cite{MW:Braess:2005, MW:JungLanger:2013}.
Discretizing problem \eqref{mw:eq:OS} by computing its weak formulations and
then inserting the finite element approximations for discretizing in space leads to the following 
discrete formulation:
\begin{align}
\label{mw:eq:OS:FE:forward}
M_h \underline{y}_{h,t} + K_h \underline{y}_h + \alpha M_h^{\Gamma_R} \underline{y}_h 
&= \beta M_h^{\Gamma_R} 
\underline{u}_h, \qquad
\underline{y}_h(0) = y_0,  \\
\label{mw:eq:OS:FE:backward}
-M_h \underline{p}_{h,t} + K_h \underline{p}_h + \alpha M_h^{\Gamma_R} \underline{p}_h &= 0, \qquad
\underline{p}_h(T) = \underline{y}_h(T) - y_d,
\end{align}
together with 
the projection formula
\begin{equation}
\label{mw:eq:OS:FE:projection}
\underline{u}_h = \mathcal{P}_{[u_a,u_b]}(-\frac{1}{\lambda} \beta 
\underline{p}_h)
\end{equation}
for $\lambda > 0$.
The problem \eqref{mw:eq:OS:FE:forward}--\eqref{mw:eq:OS:FE:projection} 
has to be solved with respect to the nodal parameter vectors
\begin{align*}
\underline{y}_h = (y_{h,i})_{i=1,\dots,n}, \quad
\underline{u}_h = (u_{h,i})_{i=1,\dots,n}, \quad
\underline{p}_h = (p_{h,i})_{i=1,\dots,n} \,\, \in \mathbb{R}^n
\end{align*}
of the finite element approximations 
$y_{h}(\boldsymbol{x}) = \sum_{i=1}^n y_{h,i} \varphi_i(\boldsymbol{x})$,
$u _{h}(\boldsymbol{x}) = \sum_{i=1}^n u_{h,i} \varphi_i(\boldsymbol{x})$
and
$p _{h}(\boldsymbol{x}) = \sum_{i=1}^n p_{h,i} \varphi_i(\boldsymbol{x})$.
The values for $\underline{y}_h$ are set to $g$ in the nodal values on the boundary $\Gamma_2$
and the problems are solved only for the degrees of freedom.
The matrices $M_h$, $M_h^{\Gamma_R}$ and $K_h$ denote the mass matrix,
the mass matrix corresponding only to the Robin boundary $\Gamma_R$
and the stiffness matrix, respectively. 
The entries of the mass and stiffness matrices are
defined by the integrals
\begin{align*}
\begin{aligned}
 M_h^{ij} = \int_{\Omega} \varphi_i \varphi_j \,d\boldsymbol{x}, \hspace{0.5cm}
 K_h^{ij} &= \int_{\Omega} \nabla \varphi_i \cdot \nabla \varphi_j \,d\boldsymbol{x}.
\end{aligned}
\end{align*}
For the time stepping we use the implicit Euler method. 
After implementing the finite element discretization and the Euler scheme, 
we apply the following projected 
gradient method:

\begin{enumerate}
\item[1.] For $k=0$, choose an initial guess $u_h^0$ satisfying the box constraints 
$u_a \leq u_h^0 \leq u_b$.
\item[2.] Solve the discrete forward problem
\eqref{mw:eq:OS:FE:forward} corresponding to \eqref{mw:eq:PDE1}--\eqref{mw:eq:PDE5} 
in order to compute $y_h^k$.
\item[3.] Solve the discrete backward problem 
\eqref{mw:eq:OS:FE:backward} corresponding to \eqref{mw:eq:PDE1red}--\eqref{mw:eq:PDE5red} 
in order to obtain $p_h^k$.
\item[4.] Evaluate the descent direction of the discrete gradient
\begin{equation}
d^k = - \nabla \bar{J}(u_h^k) = -(\beta \chi_{\Gamma_R} p_h^k + \lambda u_h^k).
\end{equation}
\item[5.] Set $u_h^{k+1} = \mathcal{P}_{[u_a,u_b]}(u_h^k + \gamma^k d^k)$ and go to step 2 unless stopping criteria
are fulfilled.
\end{enumerate}

\begin{remark}
\label{mw:remark:linesearch}
For a first implementation, the step length $\gamma = \gamma^k$ can be chosen constant for all $k$.
However, a better performance is achieved by applying a 
line search strategy as for instance
the Armijo or Wolfe conditions to obtain the best possible $\gamma^k$ in every iteration step $k$.
We refer the reader to the methods discussed 
for instance in
\cite{MW:HerzogKunisch:2010}.
However, these strategies are not subject to the present work.
\end{remark}

\section{Numerical results and conclusions}
\label{mw:sec:7:num-results}

In this section, we present 
numerical results 
for solving the type of model optimization problem discussed in this article
and draw some 
conclusions in the end. 
The numerical experiments were computed in Matlab. 
The meshes were precomputed with Matlab's \texttt{pdeModeler}.
The finite element approximation and time stepping as well as 
the projected gradient algorithm were implemented according to
the discussions in the previous two sections.

\begin{figure}[bht!] 
\includegraphics[scale=.46]{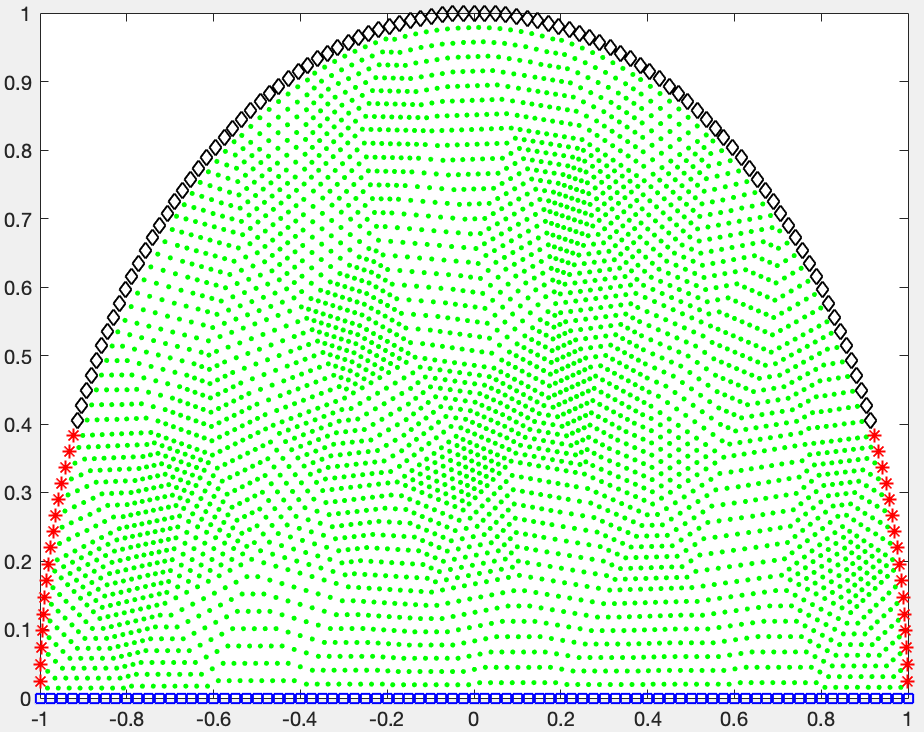}
\caption{The nodes marked corresponding to different boundaries and the interior of the domain $\Omega$ on a mesh
with 4073 nodes.}
\label{mw:fig:2:nodes}      
\end{figure}

In Figure \ref{mw:fig:2:nodes}, the nodes corresponding to the interior nodes ('g.'),
the Neumann boundary $\Gamma_1$ ('kd'), the Dirichlet boundary $\Gamma_2$ ('bs')
and the Robin boundary $\Gamma_R$ ('r*') are illustrated.

In the numerical experiments, we choose the following given data: the water temperature $g=20$,
the parameters $\alpha = \beta = 10^2$,
the final time $T = 1$,
the box constraints $u_a = 20$ and $u_b = 60$,
the desired final temperature $y_d = 30$ and the initial value
$y_0 = 0$ satisfying the boundary conditions.
For the step lengths $\gamma^k$ of the projected gradient algorithm,
we choose the golden ratio $\gamma^k = \gamma = 1.618$ 
constant for all iteration steps $k$.
The stopping criteria include that the norm of the 
errors $e_{k+1} := \|u_h^{k+1}-u_h^k\|/\|u_h^k\| < \epsilon_1$
or $|e_{k+1} - e_{k}| < \epsilon_2$ 
with $\epsilon_1 = 10^{-1}$ and $\epsilon_2 = 10^{-2}$ 
have to be fulfilled as well as setting a maximum number of iteration steps $k_{\text{max}} = 20$
with $k < k_{\text{max}}$.

In the first numerical experiment, we choose a fixed value for the cost coefficient $\lambda = 10^{-2}$
and compute the solution for different mesh sizes $n \in \{76, 275, 1045, 4073, 16081\}$.
Table \ref{mw:tab:1} presents the number of iterations needed until the stopping criteria were satisfied,
for different mesh sizes and time steps. The number of time steps was chosen
corresponding to the mesh size in order to guarantee that the CFL (Courant-Friedrichs-Lewy)
condition is fulfilled.
\begin{table}
\caption{Number of iterations needed to satisfy the stopping criteria for different mesh sizes and numbers of time steps
for a fixed cost coefficient $\lambda = 10^{-2}$.}
\label{mw:tab:1}      
\begin{tabular}{p{3cm}p{3cm}p{3cm}}
\hline\noalign{\smallskip}
mesh size & time steps & iteration steps \\
\hline\noalign{\smallskip}
 76      & 125      & 7  \\ 
 275    & 250      & 5  \\ 
1045   & 1000    & 19 \\ 
4073   & 4000    & 19 \\ 
16081 & 16000  & 4 \\ 
\hline\noalign{\smallskip}
\end{tabular}
\vspace*{-12pt}
\end{table}

\begin{figure}[bht!] 
\includegraphics[scale=.41]{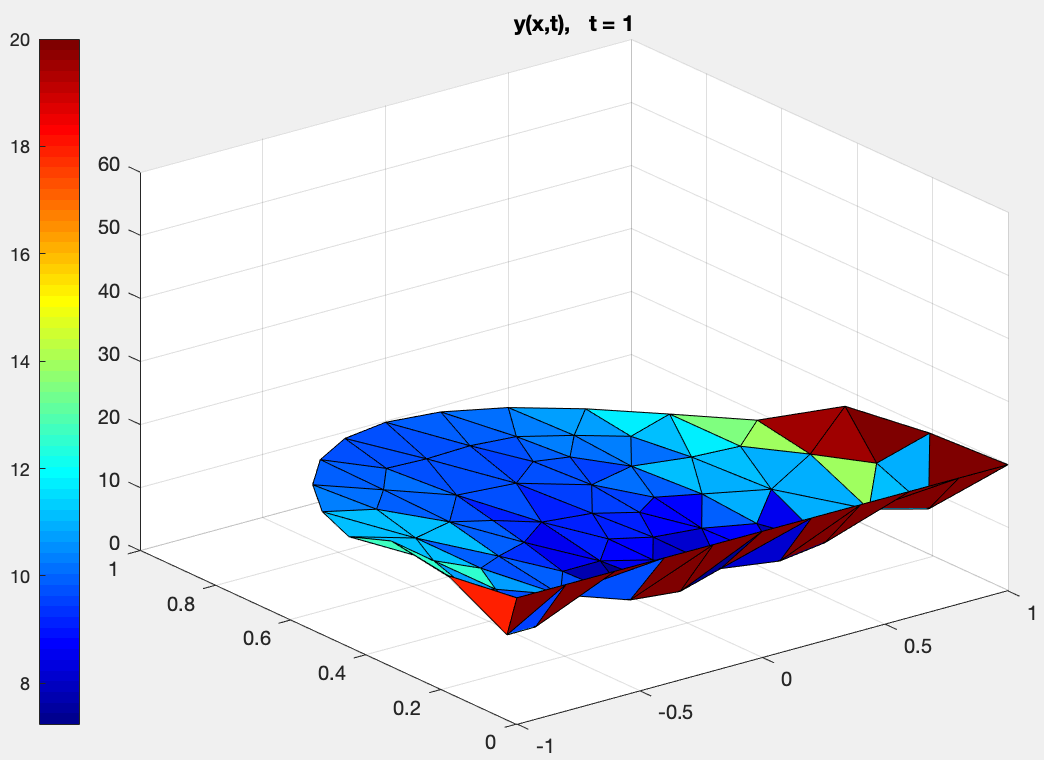}
\caption{The approximate solution $y$ for final time $t=T=1$
on a mesh
with 76 nodes.}
\label{mw:fig:3}      
\end{figure}

\begin{figure}[bht!] 
\includegraphics[scale=.41]{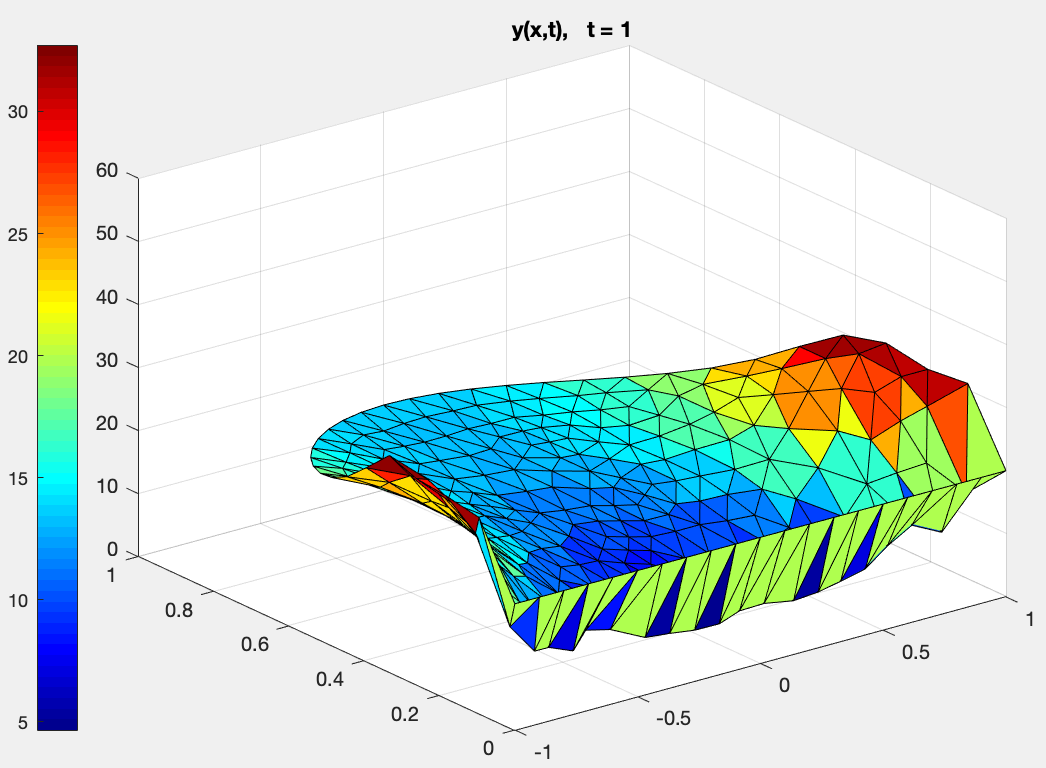}
\caption{The approximate solution $y$ for final time $t=T=1$ 
on a mesh
with 275 nodes.}
\label{mw:fig:4}      
\end{figure}

\begin{figure}[bht!] 
\includegraphics[scale=.41]{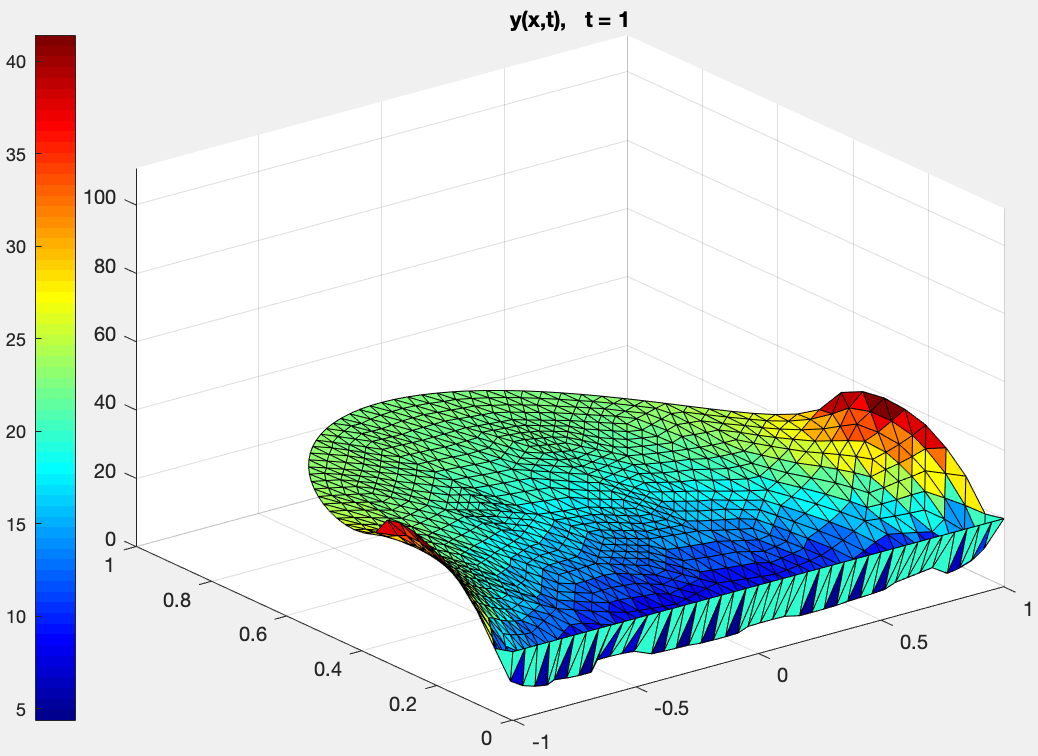}
\caption{The approximate solution $y$ for final time $t=T=1$ 
on a mesh
with 1045 nodes.}
\label{mw:fig:5}      
\end{figure}

\begin{figure}[bht!] 
\includegraphics[scale=.40]{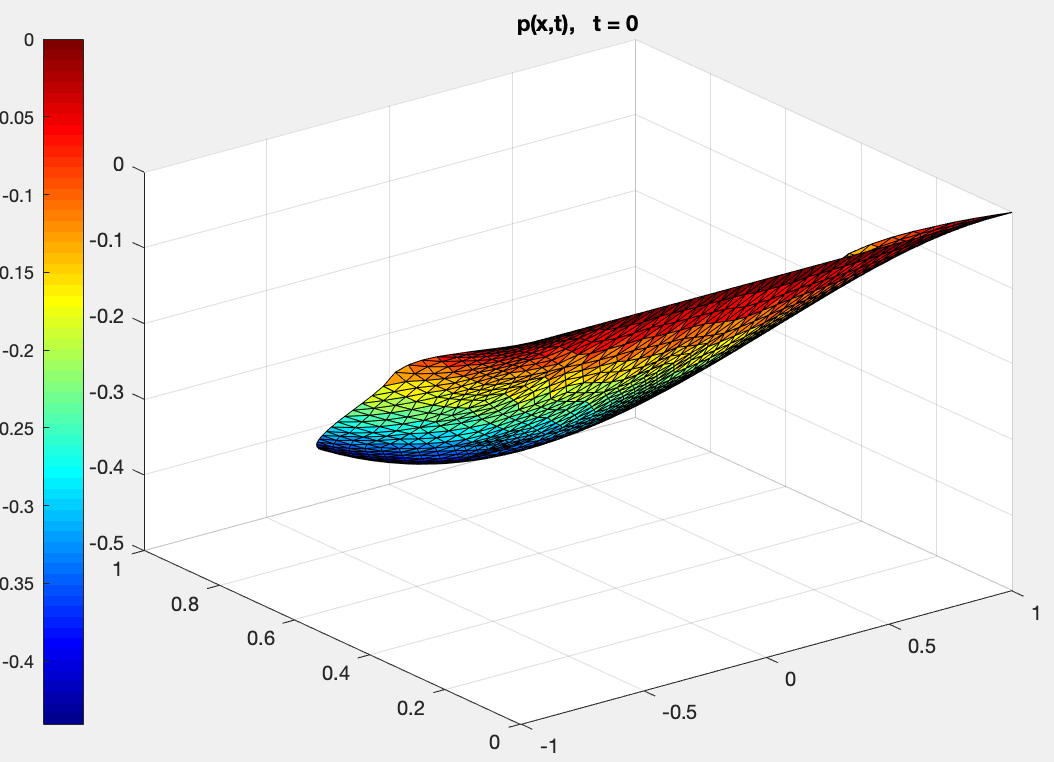}
\caption{The approximate adjoint state $p$ 
for time $t=0$ 
on a mesh with 1045 nodes.}
\label{mw:fig:8}      
\end{figure}

\begin{figure}[bht!] 
\includegraphics[scale=.41]{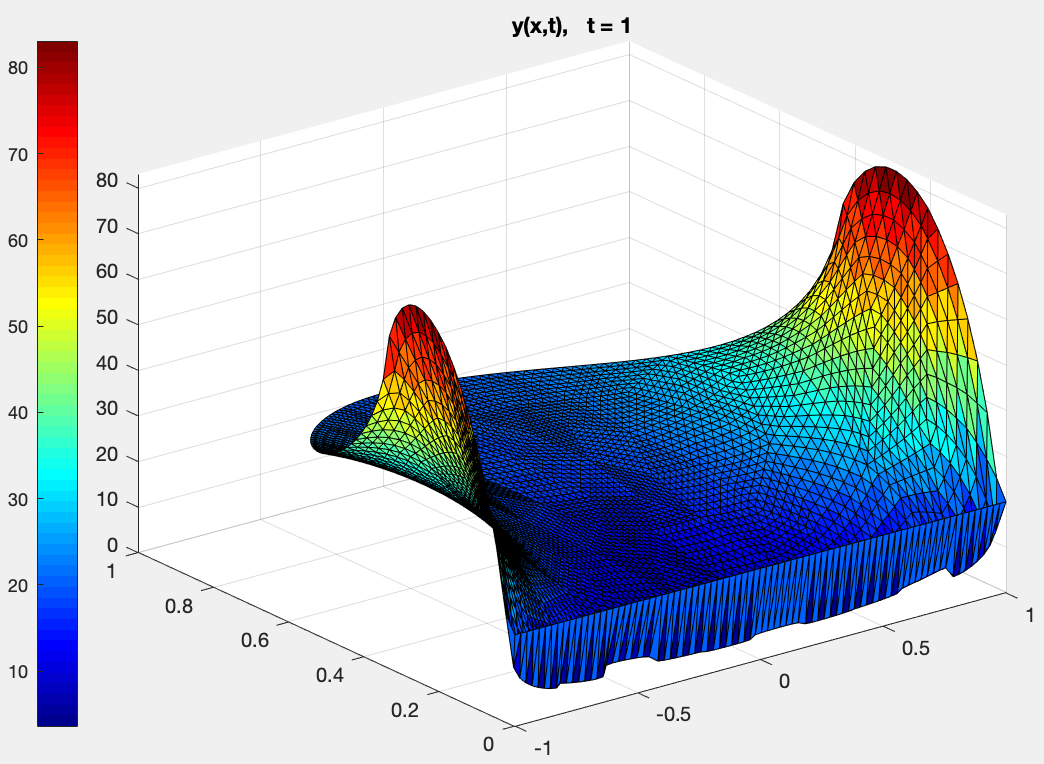}
\caption{The approximate solution $y$ for final time $t=T=1$ 
on a mesh
with 4073 nodes.}
\label{mw:fig:6}      
\end{figure}

\begin{figure}[bht!] 
\includegraphics[scale=.41]{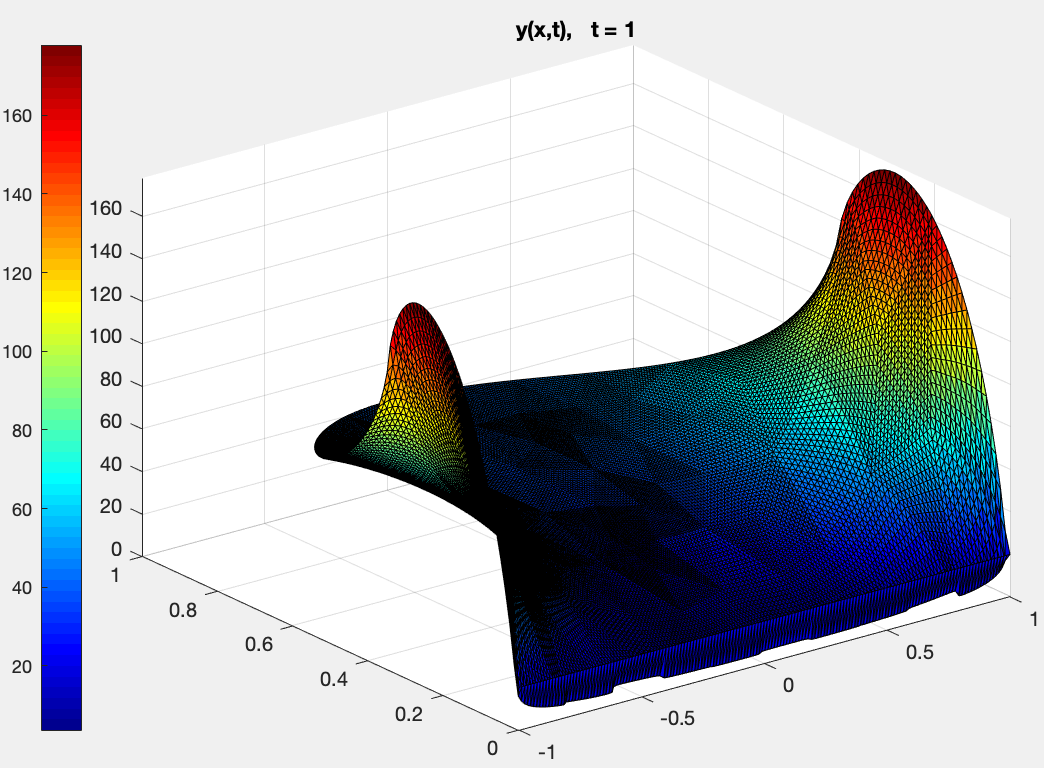}
\caption{The approximate solution $y$ for final time $t=T=1$ 
on a mesh
with 16081 nodes.}
\label{mw:fig:7}      
\end{figure}

In the set of Figures \ref{mw:fig:3} -- \ref{mw:fig:7},
the approximate solutions $y_h$ defined in Matlab as $y$ are presented for the final time $t=T=1$ computed on
the different meshes including one figure, Figure \ref{mw:fig:8}, presenting the 
adjoint state $p$ for the mesh with $1045$ nodes. We present only one figure for the adjoint state,
since for other mesh sizes the plots looked similar.
%
%
%

\begin{figure}[bht!] 
\includegraphics[scale=.41]{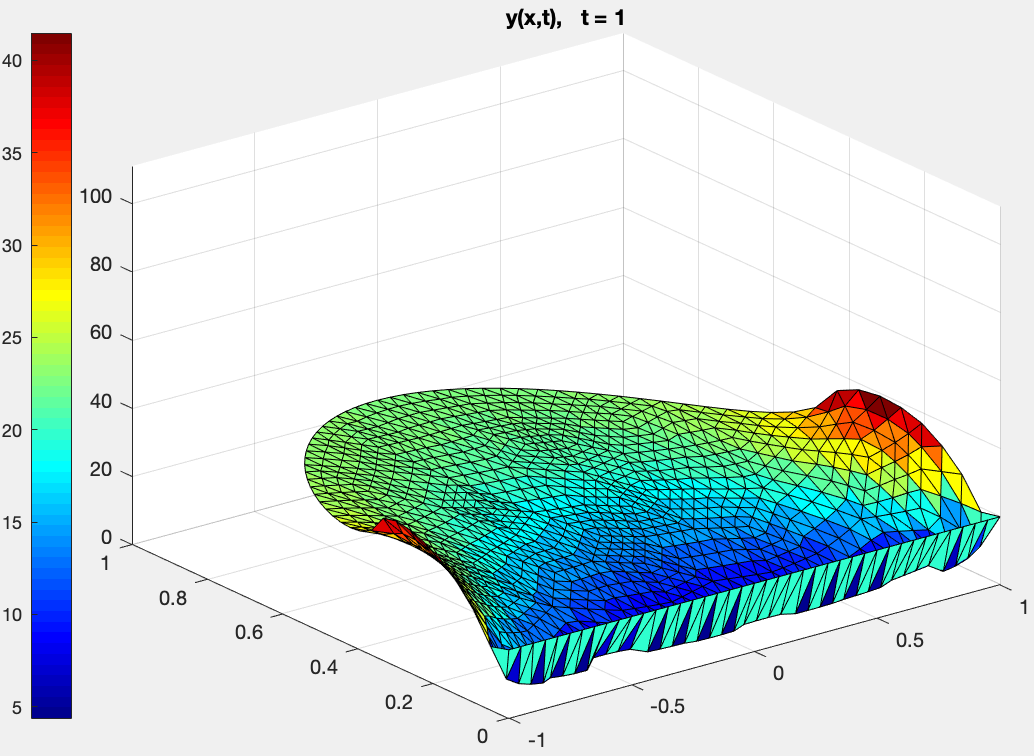}
\caption{The approximate solution $y$ for final time $t=T=1$ 
on a mesh
with 1045 nodes for the value $\lambda = 10^{-4}$.}
\label{mw:fig:9}      
\end{figure}

In the second set of numerical experiments, we compute solutions for different cost coefficients $\lambda$
on two different grids: one with mesh size $n = 275$ and $250$ time steps, and another one
with mesh size $n = 1045$ and $1000$ time steps.
The numerical results including the number of iteration steps needed are presented in Table \ref{mw:tab:2}. 
It can be observed that the numbers of iteration steps for the first case (275/250) are all very similiar
for different values of $\lambda$, even for the lower ones. 
In the second case (1045/1000), the numbers of iteration steps are getting higher the lower the values
of $\lambda$. However, the results are satisfactory for these cases too.
As example the approximate solution $y$ for the final time $t=T=1$ computed for
$\lambda = 10^{-4}$ is presented in Figure \ref{mw:fig:9}. The approximate solution
for $\lambda = 10^{-2}$ has already been presented 
in Figure \ref{mw:fig:5}. 

%
\begin{table}
\caption{Number of iterations needed to satisfy the stopping criteria for different 
cost coefficients $\lambda \in \{10^{-4}, 10^{-2}, 1, 10^{2}, 10^{4} \}$
on grids with mesh sizes $n = 275$ and $n = 1045$ with $250$ and $1000$ time steps, respectively.}
\label{mw:tab:2}      
\begin{tabular}{p{1.5cm}p{4.5cm}p{4.5cm}}
\hline\noalign{\smallskip}
$\lambda$ & iteration steps (275/250) & iteration steps (1045/1000) \\
\hline\noalign{\smallskip}
$10^{-4}$  & 5 & 19 \\ %
$10^{-2}$  & 5 & 19 \\ %
           1    & 7 & 19 \\ %
$10^{2}$   & 4 & 4 \\ %
$10^{4}$   & 4 & 4 \\ %
\hline\noalign{\smallskip}
\end{tabular}
\vspace*{-12pt}
\end{table}

The results of Tables \ref{mw:tab:1} and \ref{mw:tab:2} were included as example how one can perform
different tables for different parameter values or combinations. Students or researchers could compute exactly these
different kinds of numerical experiments in order to study the practical 
performance of the optimization
problem.

After presenting the numerical results, 
we have to mention again that the step length $\gamma = \gamma^k$ of the projected gradient
method has been chosen constant 
for all $k$ in all computations. However, better results should be achieved by applying a suitable
line search strategy, see Remark \ref{mw:remark:linesearch}.
With this we want to conclude that the optimal control problem discussed in this work
is one model formulation for solving the 
optimization of heating a domain such as the air of a swimming pool area surrounded by a glass dome.
Modeling and solving the optimal heating of a swimming pool area has potential for many different
formulations related to mathematical modeling, 
discussing different solution methods and performing many numerical tests 
including "playing around" with different values for the parameters, constants and given functions,
and finally choosing proper ones for the model problem. 
All of these tasks can be performed by students depending also on their previous knowledge, interests
and ideas.

\section*{Acknowledgements}
I would like to thank my students T. Bazlyankov, T. Briffard, G. Krzyzanowski, P.-O. Maisonneuve and C. Ne\ss ler
for their work at the 26th ECMI Modelling Week which provided part of the 
starting point 
for this work.
I gratefully acknowledge the financial support by the Academy of Finland under the grant 295897.

\end{document}